\documentclass[12pt,leqno]{article}
\usepackage{amsfonts}
\usepackage{amsmath, amsthm, amsfonts, amssymb, color}
\usepackage{mathrsfs}
\usepackage{url}
\usepackage{color}
\theoremstyle{definition}

\newcommand{\scr}[1]{\mathscr #1}
\definecolor{wco}{rgb}{0.5,0.2,0.3}

\numberwithin{equation}{section} \theoremstyle{remark}

\newcommand{\ua}{\uparrow}

\title{{\bf   Wasserstein Convergence Rate  for Empirical Measures  on Noncompact Manifolds  }\footnote{Supported in
 part by  NNSFC (11771326, 11831014, 11921001).} }
\author{
{\bf    Feng-Yu Wang$^{a),b)}$    }\\
\footnotesize{$^{a)}$ Center for Applied Mathematics, Tianjin University, Tianjin 300072, China }\\
 \footnotesize{ $^{b)}$ Department of Mathematics,
Swansea University,
Bay Campus,
Swansea,
SA1 8EN, United Kingdom}  }

\begin{document}
\allowdisplaybreaks
\def\R{\mathbb R}  \def\ff{\frac} \def\ss{\sqrt} \def\B{\mathbf
B}\def\TO{\mathbb T}
\def\I{\mathbb I_{\pp M}}\def\p<{\preceq}
\def\N{\mathbb N} \def\kk{\kappa} \def\m{{\bf m}}
\def\ee{\varepsilon}\def\ddd{D^*}
\def\dd{\delta} \def\DD{\Delta} \def\vv{\varepsilon} \def\rr{\rho}
\def\<{\langle} \def\>{\rangle} \def\GG{\Gamma} \def\gg{\gamma}
  \def\nn{\nabla} \def\pp{\partial} \def\E{\mathbb E}
\def\d{\text{\rm{d}}} \def\bb{\beta} \def\aa{\alpha} \def\D{\scr D}
  \def\si{\sigma} \def\ess{\text{\rm{ess}}}
\def\beg{\begin} \def\beq{\begin{equation}}  \def\F{\scr F}
\def\Ric{{\rm Ric}} \def\Hess{\text{\rm{Hess}}}
\def\e{\text{\rm{e}}} \def\ua{\underline a} \def\OO{\Omega}  \def\oo{\omega}
 \def\tt{\tilde}
\def\cut{\text{\rm{cut}}} \def\P{\mathbb P} \def\ifn{I_n(f^{\bigotimes n})}
\def\C{\scr C}      \def\aaa{\mathbf{r}}     \def\r{r}
\def\gap{\text{\rm{gap}}} \def\prr{\pi_{{\bf m},\varrho}}  \def\r{\mathbf r}
\def\Z{\mathbb Z} \def\vrr{\varrho} \def\ll{\lambda}
\def\L{\scr L}\def\Tt{\tt} \def\TT{\tt}\def\II{\mathbb I}
\def\i{{\rm in}}\def\Sect{{\rm Sect}}  \def\H{\mathbb H}
\def\M{\scr M}\def\Q{\mathbb Q} \def\texto{\text{o}} \def\LL{\Lambda}
\def\Rank{{\rm Rank}} \def\B{\scr B} \def\i{{\rm i}} \def\HR{\hat{\R}^d}
\def\to{\rightarrow}\def\l{\ell}\def\iint{\int}
\def\EE{\scr E}\def\Cut{{\rm Cut}}\def\W{\mathbb W}
\def\A{\scr A} \def\Lip{{\rm Lip}}\def\S{\mathbb S}
\def\BB{\scr B}\def\Ent{{\rm Ent}} \def\i{{\rm i}}\def\itparallel{{\it\parallel}}
\def\g{{\mathbf g}}\def\Sect{{\mathcal Sec}}\def\T{\mathcal T}\def\V{{\bf V}}
\def\PP{{\bf P}}\def\HL{{\bf L}}\def\Id{{\rm Id}}\def\f{{\bf f}}\def\cut{{\rm cut}}

\def\BL{\scr A}

\maketitle

\begin{abstract}  Let $X_t$ be the (reflecting) diffusion process generated by $L:=\DD+\nn V$ on a complete connected Riemannian manifold $M$ possibly with a boundary $\pp M$,
where $V\in C^1(M)$  such that  $\mu(\d x):= \e^{V(x)}\d x$ is a probability measure.  We estimate the convergence rate  for the empirical measure $\mu_t:=\ff 1 t \int_0^t \dd_{X_s}\d s$ under the  Wasserstein distance.  As a typical example, when $M=\R^d$ and $V(x)= c_1- c_2 |x|^p$ for some constants $c_1\in \R, c_2>0$ and $p>1$, the explicit upper and lower bounds are present for the convergence rate, which are of sharp order when either $d<\ff{4(p-1)}p$ or  $d\ge 4$ and $p\to\infty$.

\end{abstract} \noindent
 \vskip 2cm

\section{Introduction}

Let $M$ be a $d$-dimensional complete connected Riemannian manifold, possibly with a boundary $\pp M$. Let $V\in C^1(M)$ such that  $Z_V:= \int_M \e^{V(x)}\d s<\infty$, 
 where $\d x:={\rm vol}(\d x)$ stands for the Riemannian volume measure. Then $\mu(\d x):= Z_V^{-1} \e^{V(x)}\d x$ is a probability measure, and the (reflecting if $\pp M$ exists) diffusion process $X_t$ generated by $L:=\DD+\nn V$ is reversible with stationary distribution $\mu.$ When $M$ is compact, the convergence rate of the empirical measure
 $$\mu_t:= \ff 1 t \int_0^t \dd_{X_s}\d s,\ \ t>0$$
 under the Wasserstein distance is investigated in \cite{WZ20}. More precisely, let $\rr$ be the Riemannian distance on $M$, and  let 
$$\W_2(\mu_1,\mu_2):= \inf_{\pi\in \C(\mu_1,\mu_2)} \|\rr\|_{L^2(\pi)}$$
be the associated $L^2$-Warsserstein distance for probability measures on $M$, where $\C(\mu_1,\mu_2)$ is the class of all couplings of $\mu_1$ and $\mu_2$. 
For two positive functions $\xi,\eta$ of $t$, we denote $\xi(t)\sim \eta(t)$ if  
$c^{-1}\le \ff{\xi(t)}{\eta(t)} \le c$ holds for some constant $c>1$ and  large $t>0$.  According to \cite{WZ20},   for large $t>0$ we have 
$$ \E [\W_2(\mu_t,\mu)^2] \sim \beg{cases} t^{-1}, &\text{if}\ d\le 3,\\
t^{-1}\log t,\ &\text{if} \ d=4,\\
t^{-\ff 2 {d-2}}, \ &\text{if} \ d\ge 5,\end{cases} $$
where  the lower bound estimate on $ \E [\W_2(\mu_t,\mu)^2]$ for $d=4$ is only derived for a typical example that
$M$ is the $4$-dimensional torus and $V=0$. 
Moreover, when $\pp M$ is either convex or empty, we have 
\beq\label{CM} \lim_{t\to\infty} t \E[\W_2(\mu_t,\mu)^2] =  \sum_{i=1}^\infty \ff 2 {\ll_i^2},\end{equation} 
where $\{\ll_i\}_{i\ge 1}$ are all non-trivial eigenvalues of $-L$ (with Neumann boundary condition if $\pp M$ exists) listed in the increasing order counting multiplicities. 
See  \cite{W20, W20b} for further   studies on the conditional empirical measure of the  $L$-diffusion process with absorbing boundary. 

In this note, we investigate the convergence rate of $\E[\W_2(\mu_t,\mu)^2]$   for non-compact Riemannian manifold $M$.

\subsection{Upper bound estimate}  

We first present a  result on  the upper bound estimate  of   $\E^\nu [\W_2(\mu_t,\mu)^2]$, where $\E^\nu$ is the expectation for the diffusion process with initial distribution $\nu$. When 
$\nu=\dd_x$ is a Dirac measure, we simply denote $\E^x= \E^{\dd_x}.$

 Let $p_t(x,y)$ be the heat kernel of the (Neumann) Markov semigroup $P_t$ generated by $L$.  We will assume 
\beq\label{B1} \gg(t):= \int_M p_t(x,x) \mu(\d x)<\infty,\ \ t>0.\end{equation}
By \cite[Theorem 3.3]{W00} (see also \cite[Theorem 3.3.19]{W05}) and the spectral representation of heat kernel,  \eqref{B1} holds if and only if    $L$ has discrete spectrum such that all   eigenvalues $\{\ll_i\}_{i\ge 0}$   of $-L$ listed in the increasing order satisfy 
$$\sum_{i=0}^\infty \e^{-\ll_i t}<\infty,\ \ t>0.$$  Since $M$ is connected, the trivial eigenvalue $\ll_0=0$ is simple,  so that 
\beq\label{PI}\ll_1:= \inf\big\{\mu(|\nn f|^2):\ f\in C_b^1(M), \mu(f)=0, \mu(f^2)=1\big\}>0.  \end{equation} 
The first non-trivial eigenvalue $\ll_1$ is called the spectral gap of $L,$ and \eqref{PI} is known as the Poincar\'e inequality. 

In particular,  \eqref{B1} holds  if $P_t$ is ultracontractive, i.e.  
$$\sup_{x,y\in M} p_t(x,y)= \|P_t\|_{L^1(\mu)\to L^\infty(\mu)} <\infty,\ \ t>0.$$
Since $\gg(t)$ is decreasing in $t$, \eqref{B1} implies 
\beq\label{BB} \bb(\vv):= 1+  \int_\vv^1\d s\int_s^1 \gg(t)\d t<\infty,\ \ \vv\in (0,1].\end{equation} 
Moreover, let 
\beq\label{BD0}  \aa(\vv):=\E^\mu [\rr(X_0,X_\vv)^2]  =\int_M \rr(x,y)^2 p_\vv(x,y) \mu(\d x) \mu(\d y),\ \ \vv>0.\end{equation} 
Finally,  for any $k\ge 1$, let $\scr P_k=\{\nu\in \scr P: \nu= h_\nu\mu, \|h_\nu\|_\infty\le k\},$ where $\scr P$ is the set of all probability measures on $M$.

\beg{thm}\label{T4} Assume $\eqref{B1}$. 
\beg{enumerate} \item[$(1)$]  For any $k\ge 1$, 
\beq\label{A0} \limsup_{t\to\infty} \Big\{t \sup_{\nu\in \scr P_k} \E^\nu [\W_2(\mu_t,\mu)^2] \Big\}\le \sum_{i=1}^\infty \ff 8 {\ll_i^2}.\end{equation} 
If $P_t$ is ultracontractive, then 
\beq\label{A0'}  \limsup_{t\to\infty}  \Big\{t \E^\nu [\W_2(\mu_t,\mu)^2] \Big\}\le \sum_{i=1}^\infty \ff 8 {\ll_i^2} \end{equation} 
holds for    $\nu\in\scr P$ satisfying
\beq\label{A01}   \int_0^1  \d s\int_M  \E^\nu \big[ \rr(x,X_s)^2\big]  \mu(\d x) <\infty.\end{equation}
  \item[$(2)$]    There exists a constant $c>0$ such that    
\beq\label{B3} \sup_{\nu\in \scr P_k} \E^{\nu} \W_2(\mu_t,\mu)^2\le c k \inf_{\vv\in (0,1]}  \big\{\aa(\vv)+ t^{-1} \bb(\vv)\big\},\ \ t,k\ge 1.\end{equation}
  If $P_t$ is ultracontravtive, then there exists a constant $c>0$ such that   for any $\nu\in \scr P$ and $t\ge 1$,  
\beq\label{B3'}   \beg{split} &  \E^{\nu} [\W_2(\mu_t,\mu)^2] 
 \le c\bigg\{\ff 1 t    \int_0^1 \E^\nu \big[\mu\big(\rr(X_s,\cdot)^2\big) \big]\d s +  \inf_{\vv\in (0,1]} \big\{\aa(\vv)+ t^{-1} \bb(\vv)\big\}\bigg\}.\end{split}\end{equation} 
\end{enumerate} \end{thm}

Since the conditions \eqref{B1}, \eqref{BD0} and \eqref{A01} are less explicit, for the convenience of applications we present the following consequence of 
Theorem \ref{T4}.

\beg{cor} \label{C1} Assume that  $\pp M=\emptyset$ or $\pp M$ is convex outside a compact set. Let $V=V_1+V_2$ for some functions $V_1,V_2\in C^1(M)$ such that  
\beq\label{CVV2} \Ric_{V_1}:=\Ric-\Hess_{V_1}\ge -K,\ \ \|\nn V_2\|_\infty\le K\end{equation} holds for some constant $K>0$, where $\Ric$ is the Ricci curvature and $\Hess$ denotes the Hessian tensor.    For any $t,\vv>0$, let
\beg{align*} \tt\gg(t):= \int_M \ff{\mu(\d x)}   {\mu(B(x,\ss t))},\ \
 \tt\bb(\vv)  :=1+ \int_\vv^1\d s\int_s^1\tt\gg(r)\d r.\end{align*} 
\beg{enumerate} \item[$(1)$]  There exists a constant $c>0$ such that  
\beq\label{B3''} \sup_{\nu\in \scr P_k} \E^{\nu} [\W_2(\mu_t,\mu)^2]\le c k  \inf_{\vv\in (0,1]}  \big\{\vv + t^{-1}\tt\bb(\vv)   \big\},\ \ t,k\ge 1.\end{equation}  
\item[$(2)$] If $\|P_t \e^{\ll \rr_o^2}\|_\infty<\infty$ for $\ll,t>0$, then  for any $t\ge 1$ and $\nu\in\scr P$, 
\beq\label{B3'''}    \E^{\nu}[ \W_2(\mu_t,\mu)^2]\le c \Big[t^{-1} \nu(|\nn V|^2 )  + \inf_{\vv\in (0,1]} \big\{\vv+ t^{-1} \tt\bb(\vv)\big\}\Big].\end{equation} \end{enumerate}
\end{cor} 

\subsection{Lower bound estimate} 
Consider the modified $L^1$-Warsserstein distance
$$\tt W_1(\mu_1,\mu_2):=\inf_{\pi\in \C(\mu_1,\mu_2)}  \int_{M\times M} \{1\land \rr(x,y)\} \pi(\d x,\d y) \le \W_2(\mu_1,\mu_2).$$
We have the following result.   

\beg{thm} \label{T3}  \beg{enumerate}
\item[$(1)$] In general, there exists a constant $c>0$ such that  
\beq\label{A1-} \E^\mu[ \tt \W_1(\mu_t,\mu)^2 ] \ge c   t^{-1},\ \ t \ge 1.\end{equation}
If $\eqref{PI}$ holds, then 
 \beq\label{A1'}  \liminf_{t\to\infty} \Big\{t  \E^\nu[ \tt \W_1(\mu_t,\mu)^2] \Big\}>0,\ \ \nu\in \scr P. \end{equation}
\item[$(2)$] Let $\pp M$ be empty or convex, and let $d\ge 3$. If $\mu(|\nn V|)<\infty$ and  
\beq\label{LAA} \Ric \ge - K,\ \  V\le K\end{equation} holds for some constant $K>0$, 
then there exists a constant $c>0$ such that     
\beq\label{A2}  \inf_{\nu\in \scr P_k}  \E^\nu[ \tt W_1(\mu_t,\mu)]\ge c  (k t)^{-\ff 1 {d-2}},\ \ k,t\ge 1,\end{equation} and moreover
\beq\label{A3}   \liminf_{t\to\infty} \Big\{t^{\ff 1 {d-2}} \E^\nu [\tt W_1(\mu_t,\mu)]\Big\}>0,\ \ d\ge 4, \nu\in \scr P.\end{equation}
\item[$(3)$] Assume that $P_t$ is ultracontractive, $\pp M$ is either empty or convex, and $\Ric-\Hess_V\ge K$ for some constant $K\in\R$. Then 
\beq\label{A4}   \liminf_{t\to\infty} \inf_{\nu\in \scr P} \Big\{t^{-1} \E^\nu [W_2(\mu_t,\mu)^2]\Big\}\ge \sum_{i=1}^\infty \ff 2 {\ll_i^2}.\end{equation} 
 \end{enumerate} 
 \end{thm}

\paragraph{Remark 1.1.} According to Theorem \ref{T4}(1) and Theorem \ref{T3}(3), when $P_t$ is ultracontractive, $\pp M$ is either empty or convex, and $\Ric-\Hess_V\ge K$ for some constant $K\in\R$, we have 
$$\sum_{i=1}^\infty \ff 2 {\ll_i^2}\le \liminf_{t\to\infty}   \Big\{t^{-1} \E^\nu [W_2(\mu_t,\mu)^2]\Big\}\le \limsup_{t\to\infty}   \Big\{t^{-1} \E^\nu [W_2(\mu_t,\mu)^2]\Big\}  
\le \sum_{i=1}^\infty \ff 8 {\ll_i^2},\ \ \nu\in \scr P.$$ Because of  \eqref{CM} derived in \cite{WZ20}  in the compact setting, we may hope that the same limit formula holds for the present
non-compact setting. In particular, for the one-dimensional  Ornstein-Uhlenbeck process where $M=\R, V(x)=-\ff 1 2|x|^2$ and $\ll_i=i, i\ge 1$, we would guess 
$$\lim_{t\to\infty} \Big\{t   \E^\mu [\W_2(\mu_t,\mu)^2] \Big\} =  \sum_{i=1}^\infty \ff 2 {i^2}.$$ 
However, there is essential difficulty to prove the exact upper bound estimate as  the corresponding  calculations in \cite{WZ20} heavily depend on the estimate  $\|P_t\|_{L^1(\mu)\to L^\infty(\mu)}\le c t^{-\ff d 2}$ for some constant $c>0$ and all $t\in (0,1],$  which is available only when $M$ is compact.

\subsection{Example} 

To illustrate Corollary \ref{C1} and Theorem \ref{T3},   we consider a class of specific models, where   the convergence rate is sharp when $d<\ff{4p-1}p$ as both upper and lower bounds behave as
$t^{-1}$, and is asymptotically sharp when $d\ge 4$ and $p\to\infty$ for which  both    upper and lower bounds  are of order $t^{-\ff 2{d-2}}$. 
The assertions will be proved in Section 4.

\beg{exa}\label{Ex2}   Let $M=\R^d$ and 
 $V(x)= -\kk |x|^\aa +W(x)$ for  some constants $\kk>0, \aa>1$, and some function $W\in C^1(M)$ with $\|\nn W\|_\infty<\infty$.
\beg{enumerate}\item[$(1)$]  There exists a constant $c>0$ such that for any  $t,k\ge 1$, we have 
\beq\label{E1} \sup_{\nu\in \scr P_k} \E^{\nu} [\W_2(\mu_t,\mu)^2] \le \beg{cases} c k t^{-\ff {2(\aa-1)}{(d-2)\aa+2}}, &\text{if}\ 4(\aa-1)<d\aa,\\
c k t^{-1}\log (1+t), &\text{if}\   4(\aa-1)=d\aa,\\
 ck t^{-1}, &\text{if}\   4(\aa-1)>d\aa.\end{cases} \end{equation}
\item[$(2)$]    If $\aa>2$, then 
 there exists a constant $c>0$ such that for any $t\ge 1$, 
\beq\label{E2} \sup_{x\in \R^d} \ff{\E^{x}[ \W_2(\mu_t,\mu)^2] } {1+|x|^{2(\aa-1)}} \le    \beg{cases}c  t^{-\ff {2(\aa-1)}{(d-2)\aa+2}}, &\text{if}\ 4(\aa-1)<d\aa,\\
c t^{-1}\log (1+t), &\text{if}\   4(\aa-1)=d\aa,\\
 c t^{-1}, &\text{if}\   4(\aa-1)>d\aa.\end{cases}  \end{equation}
 \item[$(3)$]  For any probability measure $\nu$,  there exists a constant $c>0$ such that  for large $t>0$, 
 $$ \E^{\nu}[ \W_2(\mu_t,\mu)^2] \ge   \E^{\nu} [\tt \W_1(\mu_t,\mu)^2]  \ge c t^{-\ff 2 {2\lor (d-2)}}.  $$  \end{enumerate} 
  \end{exa}

\section{Proofs of Theorem \ref{T4} and Corollary \ref{C1}}
 By the  spectral representation, the heat kernel of $P_t$ is  formulated as 
\beq\label{B4} p_t(x,y)=1+ \sum_{i=1}^\infty \e^{-\ll_i t}\phi_i(x)\phi_i(y),\ \ t>0, x,y\in M,\end{equation} 
where $\{\phi_i\}_{i\ge 1}$ are the associated unit eigenfunctions with respect to the non-trivial eigenvalues $\{\ll_i\}_{i\ge 1}$ of $-L$, with the Neumann boundary condition if $\pp M$ exists. 

We will use the following inequality due to \cite[Theorem 2]{Ledoux}
\beq\label{Ledoux} \W_2(f\mu,\mu)^2\le 4 \mu(|\nn (-L)^{-1} (f-1)|^2),\ \ f\ge 0, \mu(f)=1,\end{equation} 
which is proved using an idea due to \cite{AMB}, 
see Theorem \ref{A1} below for an extension to the upper bound on $\W_p(f_1\mu, f_2\mu).$ 
To apply \eqref{Ledoux}, we   consider    the modified empirical measures 
\beq\label{B5} \mu_{\vv,t}:=f_{\vv,t}\mu,\ \ \vv>0, t>0,\end{equation}
where, according to \eqref{B4},
\beq\label{B6}  f_{\vv,t}:=\ff 1 t\int_0^t p_\vv(X_s,\cdot)= 1+\sum_{i=1}^\infty \e^{-\ll_i \vv} \xi_i(t)\phi_i,\ \ \xi_i(t):=\ff 1 t\int_0^t \phi_i(X_s)\d s.\end{equation}

\beg{proof}[Proof of Theorem \ref{T4}] (1)  It suffices to prove for   $ \sum_{i=1}^\infty \ll_i^{-2}<\infty.$
In this case, by \cite[(2.19)]{WZ20} whose proof works under the condition \eqref{B1}, we find a constant $c>0$ such that
$$\sup_{\nu\in\scr P_k} \bigg|t\E^\nu [ \mu(|(-L)^{-\ff 1 2 } (f_{\vv,t}-1)|^2)]  -\sum_{i=1}^\infty \ff 2 {\ll_i^2\e^{2\vv \ll_i}}\bigg|\le \ff{ c k }t \sum_{i=1}^\infty  \ff 1 {\ll_i^2 \e^{2\vv \ll_i}}.$$ 
This together  with \eqref{Ledoux} yields
\beq\label{*C}   t\sup_{\nu\in \scr P_k} \E^\nu [\W_2(\mu_{\vv,t},\mu)^2] 
 \le \sum_{i=1}^\infty \ff 8 {\ll_i^2} +  \ff{ c k }t \sum_{i=1}^\infty  \ff 4 {\ll_i^2},\ \ \vv>0.\end{equation}
 To approximate $\mu_t$ using $\mu_{\vv, t}$, for any $n\ge 1$ let
 $$\W_{2,n}(\mu_1,\mu_2):= \inf_{\pi\in \C(\mu_1,\mu_2)} \bigg(\int_{M\times M} \big\{n\land \rr(x,y)^2\big\} \pi(\d x,\d y)\bigg)^{\ff 1 2},\ \ \mu_1,\mu_2\in \scr P.$$
 Given $\gg\in \scr P$, let $(X_s^\gg)_{s\ge 0}$ be the (reflecting, if $\pp M\ne \emptyset$) diffusion process generated by $L$ with initial distribution $\gg$, and let $\gg P_s$ denote the distribution of $X_s^\gg$.  By the continuity of the diffusion process and the dominated convergence theorem, we have 
 $$\limsup_{\vv\downarrow 0} \W_{2,n}(\gg P_\vv, \gg)^2 =0,\ \ n\ge 1, \gg\in \scr P.$$
Observing that $\mu_{\vv,t}=\mu_tP_\vv$, we have 
 $$\limsup_{\vv\downarrow 0} \W_{2,n}(\mu_{\vv,t}, \mu_t)^2 =0,\ \ n\ge 1, t>0.$$
 Since $ \W_{2,n}(\mu_{\vv,t}, \mu_t)^2\le n$ and $\nu\le k\mu$ for $\nu\in \scr P_k$, this and the dominated convergence theorem yield 
 $$\limsup_{\vv\downarrow 0}\sup_{\nu\in \scr P_k} \E^\nu  \W_{2,n}(\mu_{\vv,t}, \mu_t)^2 \le k \limsup_{\vv\downarrow 0}  \E^\mu  \W_{2,n}(\mu_{\vv,t}, \mu_t)^2=0,\ \ n\ge 1,t>0.$$
 Combining this with \eqref{*C} and applying the triangle inequality of $\W_{2,n}$, we derive
 \beg{align*} & t\sup_{\nu\in \scr P_k} \E^\nu [\W_{2,n}(\mu_{t},\mu)^2]  \le t \limsup_{\vv\downarrow 0}\sup_{\nu\in \scr P_k} \big\{ \W_{2,n}(\mu_{\vv,t}, \mu_t)+\W_{2,n}(\mu_{t,\vv},\mu)\big\}^2 \\
 &\le \sum_{i=1}^\infty \ff 8 {\ll_i^2} +  \ff{ c k }t \sum_{i=1}^\infty  \ff 4 {\ll_i^2},\ \ \ n\ge 1, t>0.\end{align*} 
 Therefore,  for any $t>0$ we have 
\beq\label{XJ0}  t\sup_{\nu\in \scr P_k} \E^\nu [\W_2(\mu_{t},\mu)^2]=   t\sup_{n\ge 1, \nu\in \scr P_k} \E^\nu [\W_{2,n}(\mu_{t},\mu)^2] \le \sum_{i=1}^\infty \ff 8 {\ll_i^2} +  \ff{ c k }t \sum_{i=1}^\infty  \ff 4 {\ll_i^2}, \end{equation} 
which implies  \eqref{A0}. 

Next, when $P_t$ is ultracontractive, we have   
$$\dd(\vv):= \sup_{t\ge\vv, x,y\in M} p_t (x,y) <\infty,\ \ \vv>0.$$
Then the distribution $\nu_{\vv}$ of $X_\vv$ starting at $\nu$ is in the class $\scr P_{\dd(\vv)}.$ 
For any $\vv\in (0,1]$, let 
$$\bar\mu_{\vv,t}:= \ff 1 t\int_\vv^{t+\vv}\dd_{X_s}\d s.$$
 By the Markov property and  \eqref{XJ0}, we obtain  
\beq\label{XJ1}   \limsup_{t\to\infty} \Big\{t \E^\nu [\W_2(\bar \mu_{\vv,t},\mu)^2] \Big\}=  \limsup_{t\to\infty} \Big\{t\E^{\nu_{\vv} } [\W_2(\mu_{t},\mu)^2] \Big\}
\le  \sum_{i=1}^\infty \ff 8 {\ll_i^2},\ \ \vv>0. \end{equation}
On the other hand, since  
$$\pi := \ff 1 t \int_0^\vv\dd_{(X_s, X_{s+t})} \d s  + \ff 1 t\int_\vv^t \dd_{(X_s, X_s)}\d s\in \C(\mu_t, \bar\mu_{\vv,t}),  $$  
 and since the conditional distribution of $X_{s+t}$ given $X_s$ is bounded above by $\dd(1) \mu$ for $t\ge 1$, we have 
\beg{align*} & t\E^\nu[\W_2(\mu_t,\bar \mu_{\vv,t})^2 ]\le t\E^\nu\int_{M\times M} \rr(x,y)^2 \pi(\d x,\d y)\\
&= \int_0^\vv \E^\nu [\rr(X_s,X_{s+t})^2]\d s
 \le  \dd(1) \int_0^\vv \E^\nu\big[\mu\big(\rr(X_s,\cdot)^2\big)\big]\d s=: r_\vv.\end{align*} 
Combining this with \eqref{A01},  \eqref{XJ1}, and applying the triangle inequality of $\W_2$, we arrive at
\beg{align*}&\limsup_{t\to\infty} \Big\{ t \E^\nu [\W_2(\bar \mu_{t},\mu)^2]\Big\} \\
&\le\lim_{\vv\downarrow 0} \bigg( (1+r_\vv^{\ff 1 2}) \limsup_{t\to\infty}\Big\{ t \E^\nu [\W_2(\bar \mu_{\vv,t},\mu)^2] \Big\}+ (1+ r_\vv^{-\ff 1 2})r_\vv \bigg)\\
&\le  \sum_{i=1}^\infty \ff 8 {\ll_i^2}.\end{align*} 

(2) By \eqref{PI}, we have 
\beq\label{SP} \int_M |P_t f-\mu(f)|^2\d\mu\le \e^{-2\ll_1 t} \int_M|f-\mu(f)|^2\d\mu,\ \ t\ge 0, f\in L^2(\mu).\end{equation} 
By    \eqref{B4}-\eqref{B5},   and noting that $L\phi_i=-\ll_i\phi_i$ with $\{\phi_i\}_{i\ge 1}$ being orthonormal in $L^2(\mu)$, we obtain 
\beq\label{B6} \W_2(\mu_{\vv,t},\mu)^2 \le 4 \mu(|\nn (-L)^{-1}(f_{\vv,t}-1)|^2) =4\sum_{i=1}^\infty \ll_{i}^{-1} \e^{-2\ll_i\vv} |\xi_i(t)|^2.\end{equation}
Below we prove the desired assertions respectively. 

 Since for $\nu\in \scr P_k$ we have  $\E^\nu \le k \E^\mu$, it suffices to prove for $\nu=\mu$.  Since $\mu$ is $P_t$-invariant and $\mu(\phi_i^2)=1$, we have 
\beq\label{B7}   \E^{\mu}[\phi_i(X_{s_1})^2] =   \mu(\phi_i^2)=1.\end{equation} 
Next,   the Markov property yields
 $$\E^{\mu} (\phi_i(X_{s_2})| X_{s_1})= P_{s_2-s_1}\phi_i(X_{s_1})= \e^{-\ll_i(s_2-s_1)}  \phi_i(X_{s_1}),\ \ s_2>s_1.$$
Combining this with \eqref{B7} and   the definition of $\xi_i(t)$,  
 we obtain
 \beg{align*} &\E^\mu |\xi_i(t)|^2 =\ff 2 {t^2}\int_0^t\d s_1 \int_{s_1}^t \E^\mu[\phi_i(X_{s_1})\phi_i(X_{s_2})] \d s_2 \\
 &= \ff 2 {t^2}\int_0^t\d s_1 \int_{s_1}^t \E^\mu[\phi_i(X_{s_1})^2] \e^{-\ll_i(s_2-s_1)} \d s_2\le \ff{2} {t\ll_i}.\end{align*}
 Substituting into \eqref{B6} gives
\beq\label{B8}  \E^{\mu} [  \W_2(\mu_{\vv,t},\mu)^2] \le  \ff{8}t\sum_{i=1}^\infty \ll_i^{-2} \e^{-2\ll_i\vv}  
  = \ff{32}t\sum_{i=1}^\infty \int_\vv^\infty\d s \int_s^\infty  \e^{-2\ll_i r}\d r. \end{equation} 
 Noting that \eqref{SP} and the semigroup property imply
 \beg{align*} &p_{2t}(x,x)-1=\int_M |p_t(x,y)-1|^2 \mu(\d y) = \int_M |P_{\ff t 2} p_{\ff t 2} (x,\cdot)(y) -1|^2\mu(\d y) \\
 &\le \e^{-\ll_1 t} \int_M|p_{\ff t 2}(x,y)-1|^2\mu(\d y)=  \e^{-\ll_1 t} \{p_t(x,x)-1\}, \end{align*} 
 we deduce from  \eqref{B4} that 
 $$\sum_{i=1}^\infty  \e^{-2\ll_i t}=\int_{M} \big\{p_{2t}(x,x)-1\big\}\mu(\d x) \le  \e^{-\ll_1 t} \int_M \{p_t(x,x)-1\}\mu(\d x)\le \e^{-\ll_1 t}\gg(t).$$
 Therefore, by \eqref{B8} and that $\gg(t)$ is decreasing in $t$, we find a constant $c_1>0$ such that 
 \beq\label{B9} \beg{split} & \E^{\mu}   [\W_2(\mu_{\vv,t},\mu)^2 ]\le \ff{32}t \int_\vv^\infty \d s \int_s^\infty  \e^{-\ll_1 t} \gg(t)\d t\\
 &\le \ff{32}t \int_\vv^1  \bigg(\int_s^1   \gg(t)\d t +\gg(1) \int_1^\infty \e^{-\ll_ 1t}\d t\bigg) \d s+\ff{32\gg(1)}t   \int_1^\infty \d s \int_s^\infty \e^{-\ll_1 r}\d r \\
 &\le  \ff {c_1}t   \bb(\vv),\ \ \vv\in (0,1]. \end{split} \end{equation} 
 On the other hand,  \eqref{B5} and \eqref{B6}  imply  that the measure
 $$\pi(\d x,\d y):=\ff 1 t\int_0^t\big\{ \dd_{X_s}(\d x)  p_\vv(X_s,y)\mu(\d y) \big\}\d s$$
 is a coupling of $\mu_t$ and $\mu_{\vv,t}.$ Combining this   with   the fact  that $\mu$ is $P_t$-invariant, we obtain 
 \beg{align*}  \E^{\mu} [ \W_2(\mu_t,\mu_{\vv,t})^2 ]
 \le \ff {1} t\E^\mu  \int_0^t\d s\int_M \rr(X_s,y)^2p_\vv(X_s,y)\mu(\d y)   =   \aa(\vv) .\end{align*}
By \eqref{B9} and the triangle inequality of $\W_2$, this yields 
$$ \E^\mu [\W_2(\mu_t,\mu)^2]\le 2\inf_{\vv\in (0,1]} \big\{ \aa(\vv) + c_1t^{-1} \bb(\vv) \big\}.$$
Therefore,  \eqref{B3} holds for some constant $c>0$ and $\nu=\mu$. 

Finally,  let $P_t$ be ultracontractive. Then there exists a constant $c_1>0$ such that 
\beq\label{*N} \sup_{t\ge 1} p_t(x,y)\le c_1,\ \ x,y\in M.\end{equation}
 So, the distribution of $X_1$ has a distribution $\nu_1\le c_1\mu$.
Let $\bar \mu_{t} =\ff 1 t\int_0^t \dd_{X_{1+s}}\d s$. 
It is easy to see that
\beq\label{CPP} \pi := \ff 1 t \int_0^1\dd_{(X_s, X_{s+t})} \d s+ \ff 1 t\int_1^t \dd_{(X_s,X_s)} \d s\in \C(\mu_t,\bar\mu_t),\end{equation}
so that  \eqref{*N} yields
\beq\label{N1} \beg{split}  \E^\nu [\W_2(\mu_t,\bar\mu_t)^2]\le \ff 1 t \E^\nu \int_0^1 |X_s-X_{s+t}|^2\d s\le \ff {c_1}  t \E^\nu\int_0^1 \mu\big(\rr(X_s,\cdot)^2\big)\d s.\end{split}\end{equation}
On the other hand,  by   the Markov property and \eqref{B3}, we find a constant $c_2>0$ such that 
$$\E^\nu[\W_2(\bar \mu_t,\mu)^2]= \E^{\nu_1}[ \W_2(\mu_t,\mu)^2]\le c_2\inf_{\vv\in (0,1]}  \big\{\aa(\vv)+ t^{-1} \bb(\vv)\big\}.$$
Combining this with \eqref{N1} and using the triangle inequality of $\W_2$, we prove \eqref{B3'} for some constant $c>0.$ 
  \end{proof}

\beg{proof}[Proof of Corollary \ref{C1}]  (1) By  \cite[Lemma 3.5.6]{W14} and comparing $P_t$ with the semigroup generated by $\DD+\nn V_1$,  
see for instance \cite[(2.8)]{GW01}, \eqref{CVV2} implies that the Harnack inequality 
\beq\label{HI} (P_tf(x))^2\le \{P_tf^2(y) \}\e^{C+Ct^{-1}\rr(x,y)^2},\ \ x,y\in M, t\in (0,1]\end{equation} holds for some constant $C>0.$ 
Therefore, by \cite[Theorem  1.4.1]{W13} with $\Phi(r)=r^2$ and $\Psi(x,y)= C + Ct^{-1} \rr(x,y)^2$, we obtain 
$$p_{2t}(x,x)= \sup_{\mu(f^2)\le 1} (P_tf(x))^2\le \ff 1 {\int_M\e^{-C-Ct^{-1}\rr(x,y)^2}\mu(\d y) } \le \ff {\e^{3C}} {\mu(B(x,\ss{2 t}))},\ \ t\in (0,1], x\in M.$$
This implies 
\beq\label{OBS}  \gg(t)\le \e^{3C} \tt\gg(t),\ \ t\in (0,2].\end{equation} 

On the other hand, by  \eqref{CVV2} and It\^o's formula due to \cite{Kendall}, there exists constant $C_1>0$ such that 
$$\d \rr(x, X_t)^2\le \Big[C_1\big(1+\rr(x,X_t)^2 \big)+ |\nn V(x)|^2 \Big]\d t + 2\ss 2 \rr(x,X_t)\d b_t, $$
where $b_t$ is a one-dimensional Brownian motion.  So,  there exists a constant $C_2>0$ such that 
\beq\label{*D} \E^\nu [\rr(x,X_t)^2]\le (C_1+\nu(|\nn V|^2)) t\e^{C_1t}\le  C_2 (1+\nu(|\nn V|^2))t,\ \ t\in [0,1], x\in M.\end{equation} 
Then there exists a constant $c>0$ such that 
\beg{align*} & \sup_{\nu\in \scr P_k} \int_M  \E^\nu\rr(x,X_\vv)^2\mu(\d x)\le k \int_M  \E^\mu\rr(x,X_\vv)^2\mu(\d x)\\
& \le C_2k (1+\mu(|\nn V|^2))\vv \le ck\vv,\ \  \vv\in (0,1], k\ge 1.\end{align*}
  Combining this with \eqref{OBS}, we prove the first assertion by Theorem \ref{T4}(2). 
The second assertion follows from \eqref{*D} and Theorem \ref{T4}(2),  since $P_t$ is ultracontractive  provided $\|P_t \e^{\ll \rr_o^2}\|_\infty<\infty$ for $\ll,t>0$, see for instance \cite[Theorem 3.5.5]{W14}.                                                       
\end{proof}

 \section{Proof of Theorem \ref{T3}} 

  (1) We first prove that for any $0\ne f\in L^2(\mu)$,  
  \beq\label{BB0} \lim_{t\to\infty} \ff 1 t \E^\mu\bigg[\bigg|\int_0^t f(X_s) \d s\bigg|^2 \bigg]= 4\int_0^\infty \mu\big((P_s f)^2\big)\d s >0.\end{equation}
  As shown in  \cite[Lemma 2.8]{CCG} that
   the Markov property and the symmetry of $P_t$ in $L^2(\mu)$ imply 
  \beq\label{IMM} \beg{split} &\ff 1 t \E^\mu\bigg[\bigg|\int_0^t f(X_s)\d s\bigg|^2\bigg]= \ff 2 t \int_0^t \d s_1\int_{s_1}^t \E^\mu[f(X_{s_1} P_{s_2-s_1} f(X_{s_1}) ]\d s_2\\
  &=\ff 2 t \int_0^t \d s_1 \int_{s_1}^t \mu\big((P_{\ff{s_2-s_1}2} f)^2\big) \d s_2 =\ff 4 t\int_0^{t/2}   \mu\big((P_sf)^2\big)\d s\int_s^{t-s} \d r\\
  &= \ff 4 t\int_0^{t/2} (t-2s)  \mu\big((P_sf)^2\big)\d s,\ \ t>0,\end{split}\end{equation} where we have used the variable transform $(s,r)=(\ff {s_2-s_1}2, \ff{s_1+s_2}2).$ 
This implies \eqref{BB0}. On the other hand,   we take $0\ne f\in L^2(\mu)$ with $\mu(f)=0$ and $\|f\|_\infty\lor \|\nn f\|_\infty\le 1$. Then  
$$ t\E^\mu [\tt W_1(\mu_t,\mu)^2] \ge \ff 1 t   \E^\mu\bigg[\bigg|\int_0^t f(X_s)\d s\bigg|^2\bigg]. $$
Combining this with \eqref{BB0}, we prove \eqref{A1-} for some constant $c>0.$ 
  
If $\eqref{PI}$ holds, then
\beq\label{EXP} \|P_t f-\mu(f)\|_{L^2(\mu)}\le \e^{-\ll_1 t} \|f-\mu(f)\|_{L^2(\mu)},\ \ t\ge 0, f\in L^2(\mu).\end{equation}
Let $\nu=h_\nu\mu\in \scr P$ with $h_\nu\in L^2(\mu)$. Similarly to \eqref{IMM},  for any $f\in L^2(\mu)$ with $\mu(f)=0$, we have 
\beg{align*} & \ff 1 t \bigg\{\E^\nu   \bigg[\bigg|\int_0^t f(X_s) \d s\bigg|^2 \bigg]  -  \E^\mu   \bigg[\bigg|\int_0^t f(X_s) \d s\bigg|^2 \bigg] \bigg\}\\
&= \ff 1 t \int_M \{h_\nu(x)-1\}  \E^x  \bigg[\bigg|\int_0^t f(X_s) \d s\bigg|^2 \bigg]  \mu(\d x) \\
&= \ff 2 t\int_0^t \d s_1 \int_{s_1}^t \mu\big(\{h_\nu-1\} P_{s_1}\{fP_{s_2-s_1}f\}\big)\d s_2\\
&= \ff 2 t\int_0^t \d s_1 \int_{s_1}^t \mu\big(\big\{P_{s_1} (h_\nu-1)\big\}\cdot\big\{fP_{s_2-s_1}f\big\}\big)\d s_2\\ &\ge -\ff{2\|f\|_\infty} t \int_0^{t} \d s_1\int_{s_1}^t \|P_{s_1}(h_\nu-1)\|_{L^2(\mu)} \|P_{s_2-s_1} f\|_{L^2(\mu)} \d s_2.\end{align*}
Taking $0\ne f\in L^2(\mu)$ with $\mu(f)=0$ and  $\|f\|_\infty\lor\|\nn f||_\infty\le 1$, by combining this with  \eqref{BB0} and  \eqref{EXP}, we derive 
\beq\label{EX4}\beg{split} & \liminf_{t\to\infty} \Big[t \E^\nu[\tt \W_1(\mu_t,\mu)^2] \Big\}\ge  \liminf_{t\to\infty}  \bigg\{ \ff 1 t \E^\nu   \bigg[\bigg|\int_0^t f(X_s) \d s\bigg|^2 \bigg] \bigg\}\\
&\ge 4 \int_0^\infty \mu\big(|P_sf|^2\big)\d s >0,\ \  \nu=h_\nu\mu\ \text{with}\ h_\nu\in L^2(\mu).\end{split} \end{equation}
Next,   let 
$\bar\mu_{t}=\ff 1 t \int_1^{t+1}\dd_{X_s}\d s,\ t>0.$ By \eqref{CPP}    we have  
\beq\label{EXP2} \tt \W_1(\bar\mu_{t},\mu_t)\le \int_{M\times M} 1_{\{x\ne y\}} \pi(\d x,\d y) =\ff 1 {t}.\end{equation}
 Noting that for any $x\in M$ we have $\nu_{x} := p_1(x,\cdot)\mu$ with $p_1(x,\cdot)\in L^2(\mu)$,  by the Markov property and  \eqref{EX4},  we obtain
 $$ \liminf_{t\to\infty} \Big\{t \E^x [\tt \W_1(\bar \mu_t,\mu)^2] \Big\}= \liminf_{t\to\infty} \Big[t \E^{\nu_x}[\tt \W_1(\mu_t,\mu)^2] \Big\}>0.$$
 Combining this with \eqref{EXP2} and the triangle inequality  leads to
  $$ \liminf_{t\to\infty} \Big\{t \E^x [\tt \W_1(\mu_t,\mu)^2] \Big\}>0,\ \ x\in M.$$
  Therefore, by Fatou's lemma, for any $\nu\in \scr P$ we have
\beg{align*} & \liminf_{t\to\infty} \Big\{t \E^\nu [\tt \W_1(\mu_t,\mu)^2] \Big\}=  \liminf_{t\to\infty} \int_M\Big\{t \E^x [\tt \W_1(\mu_t,\mu)^2] \Big\}\nu(\d x) \\
  &\ge \int_M \Big( \liminf_{t\to\infty}  \Big\{t \E^x [\tt \W_1(\mu_t,\mu)^2] \Big\}\Big)\nu(\d x) >0,\end{align*}
which implies \eqref{A1'}.   

 (2) Let $d\ge 3$,  and let $\pp M$ be empty or convex.  
 By   \eqref{LAA},  we have  $\Ric\ge -K$  for some constant $K>0$. Then the Laplacian comparison theorem implies (see \cite{CE})
 $$ \DD \rr(x,\cdot)(y)\le \ss{K(d-1)}\coth \Big[\ss{K/(d-1)} \,\rr(x,y)\Big]\le C \rr(x,y)^{-1},\ \ (x,y)\in \hat M$$ for some constant $C>0$, 
 where $\hat M:= \{(x,y): x,y\in M, x\ne y, x\notin {\rm cut}(y)\}, $  and  ${\rm cut}(y)$ is the cut-locus of $y$.
 So, 
 $$L \rr(x,\cdot) (y)\le |\nn V(y)|+ C\big\{\rr(x,y)+\rr(x,y)^{-1} \big\},\ \ (x,y)\in \hat M.$$
 Combining this with  the It\^o's formula due to \cite{Kendall}, we obtain
 $$\d \rr(X_0,X_t)\le \ss 2 \d b_t+ \big\{|\nn V(X_t)|+ C\rr(X_0,X_t)+C\rr(X_0,X_t)^{-1} \big\}\d t +\d l_t,$$
 where $b_t$ is a one-dimensional Brownian motion, and $l_t$ is the local time of $X_t$ at the initial value $X_0$, which is an increasing process supported on $\{t\ge 0: X_t=X_0\}$. 
 Thus, we   find a constant $C_1>0$ such that 
 $$\d \Big\{\ff{\rr(X_0, X_t)^2}{1+\rr(X_0, X_t)^2} \Big\}\le C_1 (1+ |\nn V(X_t)|)\d t+\d M_t$$
 for some martingale $M_t$. Since $\mu$ is $P_t$-invariant, this implies 
 $$ \E^\mu  \big\{\rr(X_0,X_t)\land 1\big\}^2 \le C_2 \big\{1+ \mu(|\nn V|)\big\} t,\ \ t\ge 0, x\in M$$ for some constant $C_2>0$. Therefore,  for  any $N\in \mathbb N$ and $t_i:= (i-1)t/N$, 
the probability measure
  $$\tt\mu_N:= \ff 1 N \sum_{i=1}^N \dd_{X_{t_i}}=\ff 1 t\sum_{i=1}^N  \int_{t_i}^{t_{i+1}} \dd_{X_{t_i}}\d s$$ satisfies
 \beg{align*} 
 \E^\mu \tt W_1(\tt\mu_N,\mu_t)^2  &\le \ff 1 t \sum_{i=1}^N \int_{t_i}^{t_{i+1}}  \E^\mu  (\rr(X_{t_i},X_s)\land 1)^2 \d s\\
 &\le  \ff {C_3}  t \sum_{i=1}^N \int_{t_i}^{t_{i+1}} ( s-t_i) \d s  \le \ff{C_3 t} N\end{align*}   for some constant $C_3>0$. 
 So,  
 \beq\label{DD2} \sup_{\nu\in \scr P_k}  \E^\nu [\tt W_1(\tt\mu_N,\mu_t)^2]\le k  \E^\mu[ \tt W_1(\tt\mu_N,\mu_t)^2  ]  \le  \ff{C_3 kt} N,\ \ N,k\ge 1.\end{equation}     
 On the other hand, by $\Ric\ge -K$ and $V\le K$ in \eqref{LAA} and using the volume comparison theorem, 
 we find a constant $C_4>1$ such that 
 $$\mu(B(x,r))\le C_4 r^d,\ \ x\in M, r\in [0,1],$$    where 
 $B(x,r):= \{y\in M: \rr(x,y)\land 1 \le r\}$. Since $\mu$ is a probability measure, this inequality holds for all $r>0$. Therefore, by \cite[Proposition 4.2]{RE1}, 
 there exists a constant $C_5>0$ such that 
 $$\tt W_1(\tt \mu_N, \mu)\ge C_5 N^{-\ff 1 d},\ \ N\ge 1.$$
 Combining this with \eqref{DD2} and using the triangle inequality for $\tt W_1$, we obtain 
 $$\sup_{\nu\in \scr P_k}\E^\nu[ \tt W_1(\mu_t,\mu)]\ge C_5 N^{-\ff 1 d}-   \ss{C_3k t} N^{-\ff 1 2},\ \ N,k\ge 1.$$
 maximizing in $N\ge 1$, we find a constant $c>0$ such that  \eqref{A2} holds.  
 
Now, let $d\ge 4$. To prove \eqref{A3} for general probability measure $\nu$, we consider the shift empirical measure
 $$\bar\mu_t:= \ff 1 t\int_1^{t+1} \dd_{X_s}\d s,\ \ t\ge 1,$$
 and the probability measures  
  $$\nu_x:=\dd_xP_1= p_1(x,\cdot)\mu,\ \ \nu_{x,1}:= \ff{1_{B(x,1)}}{\nu_x(B(x,1))}\nu_x,\ \ x\in M.$$ By the Markov property, we obtain 
\beg{align*}& \E^x [\tt W_1(\bar\mu_t,\mu])= \E^{\nu_x} [\tt W_1(\mu_t,\mu) ]=\int_M \E^y [\tt W_1(\mu_t,\mu) ]p_1(x,y)\mu(\d y)\\
&\ge \int_{B(x,1)}  \E^y [\tt W_1(\mu_t,\mu)] p_1(x,y)\mu(\d y) =  \nu_x(B(x,1))  \E^{\nu_{x,1}}  [\tt W_1(\bar\mu_t,\mu)].\end{align*}
Noting that $h(x):= \sup_{y\in B(x,1)}p_1(x,y)<\infty$, this and \eqref{A2} yield
$$ \E^x [\tt W_1(\bar\mu_t,\mu)]\ge g(x) t^{-\ff 1 {d-2}},\ \ g(x):= c\nu_x(B(x,1)) h(x)^{-\ff 1 {d-2}},\\  x\in M, t\ge 1.$$ Consequently, for any probability measure $\nu$,
$$ \E^\nu [ \tt W_1(\bar\mu_t,\mu)]=\int_M  \E^x [\tt W_1(\bar\mu_t,\mu) ]\nu(\d x)  \ge\nu(g) t^{-\ff 1 {d-2}},\ \ t\ge 1.$$
Combining this with \eqref{EXP2}
and noting that $d\ge 4$ implies $t^{-\ff 1 {d-2}}\ge t^{-\ff 1 2}$ for $t\ge 1$, 
we find a constant $c_\nu>0$ such that when $t$ is large enough,
$$\E^\nu  [\tt W_1(\mu_t,\mu)]\ge \E^\nu \big[\tt W_1(\bar\mu_t,\mu)-\tt \W_1(\bar\mu_t,\mu_t)\big] \ge c(\nu)t^{-\ff 1 {d-2}}.$$

(3) According to \cite[Theorem 2.1]{WZ20}, for any $\vv\in (0,1]$ we have 
\beq\label{*Q1} \liminf_{t\to\infty} \Big\{t \inf_{x\in M} \E^x [\W_2(\mu_{\vv,t},\mu)^2]\Big\}\ge \sum_{i=1}^\infty \ff 2 {\ll_i^2 \e^{2\vv \ll_i}}.\end{equation} 
On the other hand, by \cite[Theorem 3.3.2]{W14}, the conditions that $\Ric-\Hess_V\ge K$ and $\pp M$ is empty or convex imply 
$$\W_2(\mu_{\vv,t},\mu)^2\le \e^{-2\vv K} \W_2(\mu_t,\mu)^2,\ \ \vv\ge 0.$$
Combining this with \eqref{*Q1}, we derive 
$$ \liminf_{t\to\infty} \Big\{t \inf_{x\in M} \E^x [\W_2(\mu_{t},\mu)^2]\Big\}\ge \e^{2\vv K} \sum_{i=1}^\infty \ff 2 {\ll_i^2 \e^{2\vv \ll_i}},\ \ \vv\in (0,1].$$
By letting $\vv\downarrow 0$ we finish the proof.

  \section{Proof of Example \ref{Ex2}}
  
    (1) Taking $V_1\in C^\infty(\R^d)$ such that $V_1(x)= -\kk |x|^\aa$ for $|x|\ge 1$, and writing $V_2= V+W-V_1$, we see that \eqref{CVV2} holds for some constant $K\in\R$.
  By Corollary \ref{C1}, it suffices to estimate $\tt\gg(t)$. For any $x\in \R^d$ with $|x|\ge 1$, and any $t\in (0,1],$   let 
  $x_t= \ff{x}{|x|}\big(|x|-\ff 1 2\ss t\big).$  We   find a constant  $c_1>0$ and some point $z\in B(x,\ss t)$ such that 
  \beq\label{BM1} \mu\big(B(x,\ss t)\big)\ge \int_{B(x_t,\ff 1 4 \ss t)} \e^{-\kk |y|^\aa+W(y)} \d y \ge c_1 t^{\ff d 2} \e^{-\kk(|x|-\ff 1 4 t^{\ff 1 2})^\aa+ W(z)}.\end{equation} 
Since $|x|\ge 1$, $t\in (0,1]$ and $\aa>1$, we find a constant $c_2>0$ such that 
\beq\label{BM2} \beg{split} &|x|^\aa -\big(|x|-t^{\ff 1 2}/4\big)^\aa = \aa\int_{|x|-\ff 1 4 t^{\ff 1 2}}^{|x|} r^{\aa-1}\d r\\
& \ge \ff{\aa t^{\ff 1 2}} 4 \Big(\ff{|x|}2\Big)^{\aa-1} \ge c_2 |x|^{\aa-1} t^{\ff 1 2}.\end{split} \end{equation}  Moreover,  
  $$|W(z)-W(x)|\le \|\nn W\|_\infty|x-z|\le \|\nn W\|_\infty,\ \ t\in (0,1], z\in B(x, t^{\ff 1 2}).$$ 
 Combining this with \eqref{BM1} and \eqref{BM2}, we find a    $c_3>0$ such that 
  $$\mu\big(B(x,\ss t)\big)\ge c_3 t^{\ff d 2}\e^{-\kk |x|^\aa + c_2 |x|^{\aa-1} t^{\ff 1 2} +W(x)},\ \ t\in [0,1], x\in \R^d.$$
Noting that $-\kk|x|^\aa+2|W(x)|$ is bounded from above,  we find  constants $c_4,c_5>0$ such that 
 $$\int_{|x|\ge 1} \ff{\mu(\d x)}{\mu(B(x,\ss t))}\le c_4 t^{-\ff d 2} \int_1^\infty r^{d-1} \e^{-c_2 r^{\aa-1} t^{\ff 1 2}}\d r\le c_5 t^{-\ff d 2 -\ff d{2(\aa-1)}}= c_5  t^{-\ff {\aa d} {2(\aa-1)} },\ \ t\in (0,1].$$
 On the other hand, there exists a constant $c_6>0$ such that $\mu(B(x, r))\ge c_6 r^{d}$ for $|x|<1$ and $r\in (0,1]$. In conclusion,  there exists a constant $c_7>0$ such that
   $$\tt\gg(t):=\int_{\R^d} \ff{\mu(\d x)}{\mu(B(x,\ss t))}\le c_5  t^{-\ff {\aa d} {2(\aa-1)} }+ c_6^{-1} t^{-\ff d 2}\le c_7 t^{-\ff {\aa d} {2(\aa-1)}},\ \  t\in (0,1].$$
Thus,  there exists a constant $c_8>0$ such that  for any $\vv\in (0,1],$ 
 $$\tt\bb(\vv)\le 1+c_6 \int_\vv^1 \d s\int_s^1t^{-\ff{d\aa}{2(\aa-1)}}\d t \le \beg{cases}c_8 \vv^{2-\ff{d\aa}{2(\aa-1)}}, &\text{if}\  2<\ff{d\aa}{2(\aa-1)},\\
c_8 \log (1+\vv^{-1}), &\text{if}\  2=\ff{d\aa}{2(\aa-1)},\\
c_8, &\text{if}\  2>\ff{d\aa}{2(\aa-1)}.\end{cases}$$
 By taking $\vv= t^{-\ff {2(\aa-1)}{(d-2)\aa+2}}$ if $4(\aa-1)<d\aa,$ $\vv= t^{-1}$ if $4(\aa-1)=d\aa,$ and $\vv\downarrow 0$ if  $4(\aa-1)>d\aa,$ 
  we derive 
\beq\label{AC0}  \inf_{\vv\in (0,1]}\big\{\vv+t^{-1}\tt\bb(\vv)\big\}\le \beg{cases} c t^{-\ff {2(\aa-1)}{(d-2)\aa+2}}, &\text{if}\ 4(\aa-1)<d\aa,\\
c t^{-1}\log (1+t), &\text{if}\   4(\aa-1)=d\aa,\\
 c t^{-1}, &\text{if}\   4(\aa-1)>d\aa\end{cases} \end{equation} for some constant $c>0$.  Therefore,   \eqref{E1} follows from  Corollary \ref{C1}(1). 
  
(2) Next,   by \cite[Corollary 3.3]{RW04},  when $\aa>2$ the Markov semigroup $P_t^0$ generated by $\DD-\kk \nn |\cdot|^\aa$ is ultracontractive
with
\beq\label{AC} \|P^0_t\|_{L^1(\mu_0)\to L^\infty(\mu_0)}\le \e^{c_1(1+t^{-\aa/(\aa-2)})},\ \ t>0 \end{equation} for some constant $c_1>0$, where $\mu_0(\d x):= Z^{-1}\e^{-\kk |x|^\aa}\d x $ is probability measure with 
normalized constant $Z>0$.  According to the correspondence between the ultracontractivity and the log-Sobolev inequality, see \cite{Davies},   \eqref{AC} holds if and only if 
there exists a constant $c_2>0$ such that 
$$\mu_0(f^2\log f^2)\le r \mu_0(|\nn f|^2) + c_2(1+ r^{-\ff \aa{\aa-2}}),\ \ r>0, \mu_0(f^2)=1.$$
Replacing $f$ by $f\e^{\ff W 2}$ and using $\|\nn W\|_\infty<\infty$ which implies $\mu(\e^{c W})<\infty$ for any $c>0$ due to $\aa>1$, we find   constants $c_3 $ such that 
\beg{align*} &\mu(f^2\log f^2)\le \mu(f^2W)+ 2r \mu(|\nn f|^2) + 2\|\nn W\|_\infty^2 +  c_2(1+ r^{-\ff \aa{\aa-2}})\\
&\le   \ff 1 2 \mu(f^2\log f^2)+ \ff 1 2 \log \mu(\e^{2W})+ 2r \mu(|\nn f|^2) +2\|\nn W\|_\infty^2+ c_2(1+ r^{-\ff \aa{\aa-2}})\\
&\le     \ff 1 2\mu(f^2\log f^2)+2r \mu(|\nn f|^2) +c_3(1+ r^{-\ff \aa{\aa-2}}),\ \ r>0, \mu(f^2)=1,\end{align*} 
where in the second line we have used the Young inequality \cite[Lemma 2.4]{ATW09} 
$$\mu(f^2 g)\le \mu(f^2\log f^2)+\log \mu(\e^g),\ \ \ \mu(f^2)=1, g\in L^1(f^2 \mu).$$
Hence, for some constant $c_4>0$ we have 
$$\mu(f^2\log f^2)\le r\mu(|\nn f|^2)+ c_4 (1+ r^{-\ff \aa{\aa-2}}),\ \ r>0, \mu(f^2)=1. $$
By the above mentioned  correspondence of the log-Sobolev inequality and semigroup estimate,  this implies  
$$ \|P_t\|_{L^1(\mu)\to L^\infty(\mu)}\le \e^{c_5(1+t^{-\aa/(\aa-2)})},\ \ t>0$$ for some constant $c_5>0$. In particular, this and $\mu(\e^{\ll |\cdot|^2} )<\infty$ 
imply $\|P_t\e^{\ll |\cdot|^2}\|_\infty<\infty$ for $t,\ll>0$, so that by Corollary \ref{C1}(2),  \eqref{E2} follows from  \eqref{AC0} and the fact that $|\nn V(x)|^2\le c'(1+|x|^{2(\aa-1)})$   holds for some constant $c'>0$.

(3)  By \cite[Corollary 1.4]{W99}, the Poincar\'e inequality \eqref{PI} holds for some constant $\ll_1>0$. Moreover, it is trivial that
the condition \eqref{LAA} holds for some constant $K\ge 0.$  So,  the desired lower bound estimate is implied by Theorem \ref{T3}. 

\paragraph{Acknowledgement.}  The author would like to thank the referees for useful comments and careful  corrections.

  \appendix
  \renewcommand{\appendixname}{Appendix~\Alph{section}}

\section{Upper\ bound\ estimate\ on \  $\W_p(f_1\mu,f_2\mu)$}
 
 \


For $p\ge 1$, let $\W_p$ be the $L^p$-Wasserstein distance induced by $\rr$, i.e.
 $$\W_p(\mu_1,\mu_2)= \inf_{\pi\in \C(\mu_1,\mu_2)} \|\rr\|_{L^p(\pi)}.$$
 According to \cite[Theorem 2]{Ledoux}, for any probability density $f$ of $\mu$, we have 
 \beq\label{APP1}  \W_p(f\mu, \mu)^p \le p^p  \mu\big( |\nn(-L)^{-1} (f-1)|^p \big).\end{equation}  The idea of the proof goes back to \cite{AMB}, in which   the following estimate is presented
 for probability density functions $f_1,f_2$:
 \beq\label{APP2} \W_2(f_1\mu_1,f_2\mu_2)^2\le \int_M \ff{|\nn (-L)^{-1} (f_2-f_1)|^2}{\scr M(f_1,f_2)}\d\mu,\end{equation} 
 where $\scr M(a,b):=1_{\{a\land b>0\}} \ff{\log a-\log b}{a-b}$ for $a\ne b$, and $\scr M(a,a)=1_{\{a>0\}}a^{-1}$.  In general,    for $p\ge 1$,   denote $\scr M_p=\scr M$ if $p=2$, and when $p\ne 2$ let 
 $$\scr M_p(a,b)=1_{\{a\land b>0\}} \ff{a^{2-p}-b^{2-p}}{(2-p)(a-b)}\text{\ for\ }a\ne b,\ \ \scr M_p(a,a)= 1_{\{a>0\}}a^{1-p}.$$ 
 In this Appendix, we   extend  estimates  \eqref{APP1} and \eqref{APP2}  as follows, which   might be useful for further studies. 
 
\beg{thm}\label{A1}  For any   probability  density functions $ f_1$ and $f_2$ with respect to $\mu$ such that  $f_1\lor f_2>0$, 
\beg{align*}  \W_p(f_1\mu, f_2\mu)^p \le \min\bigg\{ &p^p 2^{p-1} \int_M \ff{| \nn (-L)^{-1} (f_2-f_1)|^p}{(f_1+f_2)^{p-1}}\d\mu,\  p^p\int_M \ff{| \nn (-L)^{-1} (f_2-f_1)|^p}{f_1^{p-1}}\d\mu,\\
&\qquad \int_M \ff{| \nn (-L)^{-1}(f_2-f_1)|^2 }{\scr M_p(f_1,f_2)}\d\mu\bigg\}.\end{align*}
 
 \end{thm}

\beg{proof}   It suffices to prove for $p>1$. 
Let ${\rm Lip}_b(M)$ be  the set of bounded Lipschitz continuous functions  on $M$.    Consider the Hamilton-Jacobi semigroup $(Q_t)_{t>0}$ on ${\rm Lip}_b(M)$:
 $$ Q_t \phi:= \inf_{x\in M} \Big\{\phi(x)+ \ff 1 {p t^{p-1}} \rr(x,\cdot)^p\Big\},\ \ t>0, \phi\in {\rm Lip}_b(M).$$
 Then for any $\phi\in {\rm Lip}_b (M)$, $Q_0\phi:= \lim_{t\downarrow 0} Q_t\phi=\phi$, $\|\nn Q_t\phi\|_\infty$ is locally bounded in $t\ge 0$, and $Q_t \phi$ solves the  Hamilton-Jacobi equation
 \beq\label{HK0} \ff{\d}{\d t} Q_t \phi= -\ff {p-1}p  |\nn Q_t\phi|^{\ff p{p-1}},\ \ t>0.\end{equation}
Let $q=\ff{p}{p-1}.$ For any $f\in C_b^1(M)$,  and any increasing function $\theta\in C^1((0,1))$ such that $\theta_0:=\lim_{s\to 0}\theta_s=0, \theta_1:=\lim_{s\to 1}\theta_s=1$, 
by   \eqref{HK0} and the integration by parts formula, we obtain 
\beg{align*} &\mu_1(Q_1f)- \mu_2(f)= \int_0^1\Big\{\ff{\d}{\d s}  \mu\big([f_1+\theta_s(f_2-f_1)]Q_s f\big)\Big\}\d s \\
& = \int_0^1 \d s \int_M \Big\{\theta_s'(f_2-f_1) Q_s f -\ff{f_1+\theta_s(f_2-f_1)}q |\nn Q_s f|^q\Big\}\d\mu \\
&= \int_0^1\d s \int_M\Big\{\theta_s' \<\nn   (-L)^{-1}(f_2-f_1),  \nn Q_sf\> -\ff{f_1+\theta_s(f_2-f_1)}q |\nn Q_s f|^q\Big\}\d\mu \\
&\le \ff 1 p \int_M|\nn (-L)^{-1}(f_2-f_1)|^p \d\mu  \int_0^1 \ff{|\theta_s'|^p}{[f_1+\theta_s(f_2-f_1)]^{p-1}}\d s,\end{align*}
where the last step is due to   Young's inequality $ab\le a^p/p+ b^q/q$ for $a,b\ge 0$. By Kantorovich duality formula 
$$\ff 1 p\W_p(\mu_1,\mu_2)^p= \sup_{f\in C_b^1(M)} \big\{\mu_1(Q_1 f)-\mu_2(f)\big\},$$
and noting that 
\beg{align*} &f_1+\theta_s(f_2-f_1)= f_1+ f_2 - \theta_sf_1 -(1-\theta_s) f_2 \\
&= (f_1+ f_2)\Big(1-\ff{\theta_sf_1}{f_1+ f_2}- \ff{(1-\theta_s) f_2}{f_1+f_2}\Big)\\
&\ge (f_1+ f_2) \min\{1-\theta_s, \theta_s\},\end{align*}
we derive
\beq\label{ECC} \W_p(\mu_1,\mu_2)^p \le  \int_0^1\ff{|\theta_s'|^p}{\min\{\theta_s, 1-\theta_s\}^{p-1}}\d s \int_M \ff{|\nn (-L)^{-1}(f_1-f_2)|^p}{(f_1+f_2)^{p-1}}\d\mu.\end{equation}
By taking 
$$\theta_s= 1_{[0,\ff 1 2]}(s) 2^{p-1} s^p + 1_{(\ff 1 2,1]}(s) \big\{1- 2^{p-1} (1-s)^p\big\}, $$ 
which satisfies 
$$\theta_s'= p 2^{p-1}\min\{s,1-s\}^{p-1},\ \ \min\{\theta_s, 1-\theta_s\}= 2^{p-1}\min\{s,1-s\}^{p},$$ we deduce from \eqref{ECC}
that
 $$\W_p(f_1\mu, f_2\mu)^p \le p^p  2^{p-1} \int_M \ff{| (-L)^{-\ff 1 2 } (f_2-f_1)|^p}{(f_1+f_2)^{p-1}}\d\mu. $$
Next, \eqref{ECC} with $\theta_s= 1-(1-s)^p$ implies 
$$\W_p(f_1\mu, f_2\mu)^p \le p^p \int_M \ff{| (-L)^{-\ff 1 2 } (f_2-f_1)|^p}{f_1^{p-1}}\d\mu.$$
Finally, with $\theta_s=s$ we deduce from \eqref{ECC}  that
$$\W_p(f_1\mu, f_2\mu)^p \le \int_M \ff{|(-L)^{-\ff 1 2}(f_2-f_1)|^2 }{ \scr M_p(f_1,f_2)}\d\mu.$$
Then the proof is finished. 
\end{proof}

\begin{thebibliography}{999}
      
  \bibitem{ATW09} M. Arnaudon, A. Thalmaier, F.-Y. Wang,  \emph{Gradient estimates and Harnack inequalities on non-compact Riemannian manifolds, }  Stoch. Proc. Appl. 119(2009), 3653--3670.
 \bibitem{AMB}  L. Ambrosio, F. Stra, D. Trevisan, \emph{A PDE approach to a 2-dimensional matching problem,} Probab. Theory Relat. Fields
  173(2019), 433--477.
  



\bibitem{CCG} P. Cattiaux, D. Chafa{\i}, A. Guillin,  \emph{ Central limit
theorems for additive functionals of ergodic Markov diffusions
processes,}    ALEA Lat. Am. J. Probab. Math. Stat. 9(2012), 337--382.
\bibitem{CE}  J. Cheeger and D. G. Ebin, \emph{Comparison Theorems in
Riemannian Geometry,} Amsterdam: North-Holland, 1975.
\bibitem{Davies} E. B. Davies,  \emph{Heat Kernels and Spectral Theory,} Cambridge University Press, 1989.

\bibitem{GW01}  F. Z. Gong,  F.-Y. Wang,  \emph{Heat kernel estimates with applications to compactness of manifolds,}  Quart J. Math. 52(2001), 1--10.
\bibitem{Kendall} W. S. Kendall, \emph{The radial part of Brownian motion on a manifold: a semimartingale property,} Ann. Probab. 15(1987), 1491--1500. 
\bibitem{RE1} B. Kloeckner, \emph{Approximation by finitely supported measures,}  ESAIM Control Optim. Calc. Var. 18(2012), 343--359.
\bibitem{Ledoux} M. Ledoux,  \emph{On optimal matching of Gaussian samples},  Zap. Nauchn. Sem. S.-Peterburg. Otdel. Mat. Inst. Steklov. (POMI) 457, Veroyatnost' i Statistika.  25, 226--264 (2017).


\bibitem{RW04}  M. R\"ockner,  F.-Y. Wang, \emph{Supercontractivity and ultracontractivity for (non-symmetric) diffusion semigroups on manifolds, } Forum Math. 15(2003), 893--921. 



\bibitem{W99} F.-Y. Wang, \emph{Existence of the spectral gap for elliptic operators, }  Arkiv f\"or Math. 37(1999), 395--407.
 
\bibitem{W00}  	F.-Y. Wang,  \emph{Functional inequalities, semigroup properties and spectrum estimates,}  Infinite Dimensional Analysis, Quantum Probability and Related Topics 3:2(2000), 263--295.
\bibitem{W10} F.-Y. Wang,  \emph{Harnack inequalities on manifolds with boundary and applications,}  J. Math. Pures Appl. 94(2010), 304--321. 


\bibitem{W05} F.-Y. Wang, \emph{Functional Inequalities, Markov Semigroups and Spectral Theory,}  2005 Science Press. 
\bibitem{W13} F.-Y. Wang, \emph{ Harnack Inequality  for Stochastic Partial Differential Equations, } Math. Brief. Springer, 2013
\bibitem{W14} F.-Y.  Wang, \emph{Analysis for Diffusion Processes on Riemnnian Manifolds, } World Sicentific, 2014.


\bibitem{W20}  F.-Y. Wang, \emph{Precise limit in  Wasserstein distance for conditional empirical measures of Dirichlet diffusion processes,} J. Funct. Anal. 280(2021), 108998.

\bibitem{W20b}  F.-Y. Wang, \emph{Convergence in Wasserstein distance for empirical measures of Dirichlet diffusion processes on manifolds,} arXiv:2005.09290.
  
 \bibitem{WZ20} F.-Y. Wang, J.-X. Zhu,  \emph{Limit theorems in Wasserstein distance for empirical measures of diffusion processes on Riemannian manifolds,} aXiv:1906.03422.


\



















 











 \end{thebibliography}
\end{document}